\documentclass[journal]{IEEEtran}

\ifCLASSINFOpdf
\else
\fi
\begin{document}
%
\title{Improved Hoeffding's Lemma and   Hoeffding's  Tail Bounds}
%
%
%

\author{David~Hertz,~\IEEEmembership{Senior Member,~IEEE}
\thanks{D. Hertz is retired.
 e-mail: hertz@013.net.il or drdavidhertz@gmail.com.}
\thanks{Manuscript  submitted on October  3, 2020 to IEEE Signal Processing Letters.}}

%
%

\markboth{IEEE Signal Processing Lettrs,~Vol.~?, No.~?, Month~Year}%
{Shell \MakeLowercase{\textit{et al.}}: Bare Demo of IEEEtran.cls for Journals}
%



\maketitle

\begin{abstract}
The purpose of this letter is   to   improve Hoeffding's lemma and consequently  Hoeffding's tail bounds. The improvement pertains  to left skewed zero mean random variables
$X\in[a,b]$, where $a<0$ and $-a>b$.
The proof of   Hoeffding's improved lemma uses Taylor's expansion, the convexity of $\exp(sx), s\in {\bf R}$ and an unnoticed  observation since Hoeffding's  publication in 
1963 that for $-a>b$ the maximum of  the intermediate function $\tau(1-\tau)$ appearing in Hoeffding's proof is attained 
 at an endpoint rather than at $\tau=0.5$ as in the case $b>-a$.
Using Hoeffding's  improved lemma we obtain one sided and two sided  tail bounds  for $P(S_n\ge t)$ and $P(|S_n|\ge t)$, respectively, where 
 $S_n=\sum_{i=1}^nX_i$ and the  $X_i\in[a_i,b_i],i=1,...,n$ are  independent zero mean  random variables (not necessarily identically distributed). 
It is interesting to note that we  could  also improve Hoeffding's two sided bound for all $\{X_i:  a_i\ne b_i,i=1,...,n\}$. This is  so because here the
one sided bound should be  increased by 
$P(-S_n\ge t)$,  wherein the left skewed intervals become right skewed and vice versa.

\end{abstract}


%
\IEEEpeerreviewmaketitle

\section{Introduction}
%
%
%
%
\IEEEPARstart{B}{y} googling Hoeffding and Signal Processing we obtain many applications and theoretical uses of Hoeffding's bounds to signal processing, e.g., to time series analysis, compressed sensing, sensor networks, financial signal processing, to name but  a few. Other applications are to  machine learning, communication and information theory \cite{bib:Igal}, randomized algorithms \cite{bib:MU}, to name but  a few. Communication networks are treated in  \cite{bib:MU}  where  the authors re derive Hoeffding's results   and present applications to packet routing in sparse networks. Our starting point is to present Hoeffding's lemma without proof somewhat differently than the original one  and then present the improved Hoeffding's lemma and prove it.
Most of the proof is not new and presented here for the sake of completeness. When the underlying  distribution is skewed to the left we 
show how the bound in Hoeffding's Lemma can be improved.

We assume that  $X$ is a  zero mean real valued random variable such that $X\in[a,b], a<0, b>0$.
Then, Hoeffding's lemma \cite{bib:WH} states that for all  $s\in \bf R, s>0$,

\begin{equation}
\label{eqn:H}
E[e^{sX}]\le \exp\left(\frac{s^2(b-a)^2}{8}\right).
\end{equation}
Let $A$ and $G$ denote the arithmetic and geometric means of $|a|$ and $b$, respectively. Noting that $b-a=b+|a|$ we have

\begin{equation}
\label{eqn:A}
A:=\frac{b+|a|}{2}
\end{equation}
and

\begin{equation}
\label{eqn:G}
G:=\sqrt{|a|b}.
\end{equation}

Hence,
 \begin{equation}
\label{eqn:equiv}
 \exp\left(\frac{s^2 (b-a)^2}{8}\right)\equiv \exp\left(\frac{s^2 A^2}{2}\right). 
\end{equation}
The improvement in  Hoeffding's bound  pertains to the case $|a|>b$. We have
	
\begin{equation}
\label{eqn:H}
E[e^{s X}]\le
\left\{
    \begin{array}{rl}
        \exp\left(\frac{s^2 A^2}{2}\right), &   b\ge |a|, \\
        \exp\left(\frac{s^2 G^2}{2}\right), &  b\le |a|,
    \end{array}
\right.
\end{equation}
where since $A\ge G$ gives for $-a>b$  a tighter bound than for $-a<b$. The proof will be given in the next Section.
The organization of the remaining Sections  is as follows. in Section~\ref{sec:PIH} we present the proof of  Hoeffding's improved lemma.  In 
Section~\ref{sec:HB} we present   Hoeffding's improved one sided tail bound and its proof. In 
Section~\ref{sec:absHB} we present   Hoeffding's improved two sided tail bound and its proof.
Finally, in  Section~\ref{sec:fin} we give the conclusion.

\section{Proof of   Hoeffding's  Improved Lemma}
\label{sec:PIH}
Since $e^{sx}$ is a convex function of $x$ and  $s>0$ is a parameter we obtain
\begin{equation}
e^{sx}\le \frac{b-x}{b-a}e^{sa}+\frac{x-a}{b-a}e^{sb}
\end{equation}
Let 
\begin{equation}
\lambda:=\frac{-a}{b-a}
\end{equation}
and
\begin{equation}
\label{eqn:u}
u:=s(b-a)>0.
\end{equation}

Hence, since ${\bf E}X=0$, using and the convexity of $\exp(sx)$ and some algebra  we obtain

\begin{equation}
\label{eqn:cher}
\begin{array}{ll}
{\bf E}e^{sX} &\le  \frac{b-{\bf E}X}{b-a}e^{sa}+ \frac{{\bf E}X-a}{b-a}e^{sb} \\         
                   & =    (1-\lambda)e^{sa}+\lambda e^{sb}\\
                   & =    (1-\lambda+\lambda e^{u}) e^{-\lambda u}\\
                   & =    e^{\psi(u)},

	\end{array}
\end{equation}

where
\begin{equation}
\psi(u):=-\lambda  u+\ln(1-\lambda+\lambda e^u).
\end{equation}
Since $u>0$ and $\psi(u)$ is well defined, by Taylor's expansion we obtain
\begin{equation}
\psi(u)=\psi(0)+\psi'(0)u+0.5\psi''(\mu)u^2, {\rm ~for~some~} \mu\in[0,u].
\end{equation}
We have
\begin{equation}
\psi'(u)=-\lambda+\frac{\lambda e^u}{1-\lambda+\lambda e^u}
\end{equation}
and
\begin{equation}
\psi''(\mu)=\frac{\lambda e^\mu}{1-\lambda+\lambda e^\mu}\left(1-\frac{\lambda e^\mu}{1-\lambda+\lambda e^\mu}\right).
\end{equation}
Let
\begin{equation}
\label{eqn:tau}
\tau(\mu)=\frac{\lambda}{(1-\lambda) e^{-\mu}+\lambda},  \mu\in[0, u].
\end{equation}

Since
$\psi(0)=\psi'(0)=0$ we obtain
\begin{equation}
\psi(u)=0.5\tau(1-\tau)u^2.
\end{equation}

Using (\ref{eqn:tau}) and the fact $u>0$ we obtain that  $\tau\in[\lambda,1]$, where the endpoints $\lambda$  and $1$ correspond to  $u=0$ and $u=\infty$, respectively. 
Now, for  $-a\le b$ since $\lambda\le 0.5$  the maximum of $\tau(\tau-1)$ is attained at $\tau=0.5$ 
and is $0.25$ giving  Hoeffding's original lemma. However, for  $-a>b$ , $\lambda> 0.5$   the maximum is attained at $\tau=\lambda$. 
Using (\ref{eqn:u}) we obtain after some algebra
\begin{equation}
\psi(u) \le 0.5\lambda(1-\lambda)u^2 =0.5s^2G^2.
\end{equation}
Hence, using (\ref{eqn:cher}) we obtain
\begin{equation}
E[e^{s X}]\le    \exp\left(\frac{s^2 G^2}{2}\right), -a\ge b.
\end{equation}

This completes the proof.

\section{   Hoeffding's Improved One Sided Tail Bound}
\label{sec:HB}

Let  $X_1,...,X_n$ be   independent random variables such that  $X_i\in[a_i,b_i], a_i<0,b_i>0~{\rm and}~{\bf E}X_i=0~{\rm for}~ i =1,...n$.
Let $S_n:=\sum_{i=1}^nX_i$ then ${\bf E}S_n=0$ . For all $s>0$ we have
\begin{equation}
\label{eqn:HH}
\begin{array}{lll}
  P(S_n\ge t)&=P\left(e^{sS_n}\ge e^{st}\right ) &{\rm Chernoff}\\
                           &\le  e^{-st}{\bf E} e^{sS_n} &{\rm  Markov's~ inequality} \\
                           &=  e^{-st}\Pi_{i=1}^n {\bf E}e^{sX_i}. & X_i{\rm 's~are~independent}.
\end{array}
\end{equation}
Now define
\begin{equation}
\label{eqn:I}
I:=\{i: b_i\ge-a_i, i=1:n\} 
\end{equation}
and
\begin{equation}
\label{eqn:J}
J:=\{i: b_i<-a_i, i=1:n\}. 
\end{equation}
Hence, using (\ref{eqn:HH})   we obtain

\begin{equation}
\begin{array}{ll}
  P(S_n\ge t)   &\le  e^{-st}\Pi_{i=1}^n {\bf E}e^{sX_i}\\
                                    &=\exp\left(-st+0.5s^2\left(\sum_{i\in I}A_i^2+\sum_{j\in J}G_j ^2\right)\right) \\
                                    &=\exp(-st+0.5s^2M_n^2),
\end{array}
\end{equation}

where  the mixed sum $M_n^2$ is defined by
\begin{equation}
M_n^2:=\sum_{i\in I}A_i^2+\sum_{j\in J}G_j^2.
\end{equation}
Finally, minimizing the exponent in the last inequality we obtain
\begin{equation}
  P(S_n\ge t)   \le \exp(-0.5(t/M_n)^2).
\end{equation}
Now let 
\begin{equation}
\label{eqn:bM}
\bar M_n^2:=\frac{M_n^2}{n},
\end{equation}
\begin{equation}
\bar A_n^2:=\frac{\sum_{i=1}^nA_i^2}{n},
\end{equation}
and
\begin{equation}
\bar G_n^2:=\frac{\sum_{i=1}^nG_i^2}{n}.
\end{equation}
Since  $G_i \le A_i, i=1,...n$ we obtain
\begin{equation}
\bar G_n\le\bar M_n\le\bar A_n.
\end{equation}
Consequently,  we obtain the new bound
\begin{equation}
 \label{eqn:HB}
\begin{array}{lll} P(S_n\ge t) &  \le  \exp(-0.5(t/\bar M_n)^2n) & {\rm Improved~bound}\\
                                                    &     \le  \exp(-0.5(t/\bar A_n)^2n)& {\rm Original~bound}.
\end{array}
\end{equation}
If we define $k$ by  $k:=t/\bar M_n$ we obtain
\begin{equation}
\begin{array}{ll}
 P(S_n/\bar M_n\ge k)&=P\left(\sum_{i=1}^n(X_i/\bar M_n)\ge k\right) \\
                              &\le  \exp(-0.5k^2n) .
\end{array}
\end{equation}
Hence, $\bar M_n$ is a  natural unit for  $X_i, i=1,...,n$ and $S_n$ that we  denote  by  $[X_i]=[S_n]=\bar M_n$.  Using this unit  we obtain

\begin{equation}
 P(S_n\ge k) \le  \exp(-0.5k^2n) .
\end{equation}

Further, note that if  $X_i \in [a,b], i=1,...,n$ are zero mean independent random  variables (not necessarily identically distributed) then using (\ref{eqn:A}) and (\ref{eqn:G}) we obtain
\begin{equation}
\label{eqn:idHH1}
  P(S_n\ge t) \le     \exp\left(-0.5(t/A)^2n\right), {\rm ~if~ }-a\le b 
\end{equation}
and
\begin{equation}
\label{eqn:idHH2}
  P(S_n\ge t) \le    \exp\left(-0.5(t/G)^2n\right), {\rm ~if~ }-a>b.
\end{equation}

\section{ Hoeffding's Improved Two Sided Bound }
\label{sec:absHB}

We will prove that
\begin{equation}
\label{eqn:idnHH}
P(|S_n|\ge t) \le
 \exp(-0.5(t/\bar M_n)^2n)+\exp(-0.5(t/\bar N_n)^2n),
\end{equation}

where $\bar M_n^2$ is as defined in (\ref{eqn:bM}) and $\bar N_n^2$ is a sort of its complement. I.e., 
\begin{equation}
\bar N_n^2:=N_n^2/n,
\end{equation}
where
\begin{equation}
N_n^2:=\sum_{i\in J}A_i^2+\sum_{j\in I}G_j^2,
\end{equation}
and $I, J$  are as defined in (\ref{eqn:I}) and (\ref{eqn:J}), respectively.

{\bf Proof of (\ref{eqn:idnHH}).}\\
Since $S_n=\sum_{i=1}^n X_i, X_i\in[a_i,b_i]$ it follows that $-S_n=\sum_{i=1}^n (-X_i),  (-X_i)\in[-b_i,-a_i], -b_i<0, -a_i>0$. Hence,
 in this case we should  use   $\bar N_n$  instead of $\bar M_n$ and  we obtain
\[
\begin{array}{ll}
P(|S_n|\ge t) &=P(S_n\ge t {\rm~or~} -S_n\le t) \\
                    &=P(S_n\ge t)+P(-S_n\ge t)       \\       
                    &\le  \exp(-0.5(t/\bar M_n)^2n) +\exp(-0.5(t/\bar N_n)^2n).
\end{array}
\]
This completes the proof. \\ Note that unless $\bar M_n=\bar N_n$ we can not define $k$ as in the one sided case.
It is easily seen and interesting to  note that if $[a_i,b_i]=[a,b], i=1,...n$,     we obtain
\begin{equation}
P(|S_n|\ge t) \le \exp(-0.5(t/G)^2n) +  \exp(-0.5(t/A)^2n).
\end{equation}
Hence, unless $A=G$ we thus  improved Hoeffding's two sided tail bound.

\section{Conclusion}
\label{sec:fin}
In this letter we presented  Hoeffding's improved lemma and Hoeffding's  improved one and two sided tail bounds for bounded random variables. The improvement pertains only to
intervals $[a,b]$, where $-a>b$ and are associated with the geometric mean of $-a$ and $b$.  For  $-a<b$  Hoeffding's lemma and Hoeffding one sided tail bound remain intact and   are associated only with the average of $-a$ and $b$.
It is interesting to note that we  could  also improve Hoeffding's two sided bound for all $\{X_i:  a_i\ne b_i,i=1,...,n\}$. This is  so because here the
one sided bound should be  increased by the bound for
$P(-S_n\ge t)$,  wherein  left skewed intervals become right skewed and vice versa.
Perhaps,   further research will focus on trying to improve other inequalities that use   Hoeffding's results.

\end{document}